\documentclass[11pt]{article}
\usepackage{graphics,amssymb,amsmath,epsf}

\usepackage{fancyhdr}

\begin{document}

\date{}
\title{On confidence intervals for the power of F-tests}
\author{A. A. Jafari*,  A. Bazargan-Lari$^\dag$ , Mingfei
Li$^\ddag$\\{\small *Department of Statistics, Shiraz University,
Shiraz 71454, IRAN}
\\ $^\dag${\small Department of Statistics, Islamic Azad University-Fars Science \& Research Branch, Shiraz,
IRAN}
\\ {\it E-mail: Bazargan-lari@susc.ac.ir}
\\ $^\ddag${\small Department of Mathematics, Michigan State University,
East Lansing, MI 48824}
\\ \it{E-mail:  limingfe@msu.edu}}
\maketitle

\begin{abstract}
This note points out how confidence interval estimates for
standard deviation transform into confidence interval estimates
for the power of F-tests at fixed alternative means. An
application is shown for the test of a two-sided hypothesis for
the mean of a normal distribution.
\end{abstract}

\bigskip


\bigskip

Let $Y$ be a random variable distributed $P_\sigma $ and
$(a(Y),b(Y))$ be a $100(1 -\gamma )\%$ confidence interval
estimator for the parameter $\sigma $. If $\omega=f(\sigma )$ is a
strictly increasing function of $\sigma $, then
$(f(a(Y)),f(b(Y)))$ is a $100(1 -\gamma )\%$ confidence interval
estimator for $\omega $ and if $\omega= f(\sigma )$ is a strictly
decreasing function of $\sigma $, then $(f(b(Y)),f(a(Y)))$ is a
100(1 -- \textit{$\gamma $})\% confidence interval estimator for
$\omega $.

Tarasin\'{s}ka (2005) uses the above idea in discussing
confidence interval estimation of the power of the one-sided
t-test based on confidence interval estimation of standard
deviation. We point out how the power of the F-test is easily
handled as well, a result that subsumes the two-sided t-test.

Consider the power $\omega $ of the typical$\,$ $\alpha $ level
F-test of the mean vector based on $X_1, X_2, \dots , X_n$
independent, normally distributed random variables with common
standard deviation $\sigma $. Under the standard model, the
F-statistic has the distribution of $(U/u)/(V/v)$ where $U$ is
distributed noncentral chi-square with $u$ degrees of freedom and
noncentrality parameter $\delta $, $V$ is distributed central
chi-square with $v$ degrees of freedom, and $U$ and $V$ are
independent.  At fixed alternative mean vectors, the
noncentrality parameter $\delta=\delta (\sigma)$ is a strictly
decreasing function of the scale parameter $\sigma $. (See
Scheff\'{e} (1959, p. 39).)  As an example consider the two-way,
normal theory ANOVA with $K$ observations in each of $IJ$ cells,
and the F-test that all interactions are 0.  Here $u =(I-1)(J-1)$
, $v=(K- 1)IJ$ and the noncentrality parameter is
$\delta(\sigma)=\sqrt{K\sum (\alpha \beta )_{ij}^{2}  } /\sigma $.

Let $G_{u,v,\delta } $ denote the cdf of $(U/u)/(V/v)$, that is,
the noncentral F-distribution.  It is easy to see from the
representation $(U/u)/(V/v)$ that
\begin{equation}
Power=\omega=1-G_{u,v,\delta } (c)=1-P(U\le
cuV/v)=1-E[F_{u,\delta } (cuV/v)]
\end{equation}
where the expectation $E$ is with respect to $V$, $F_{u, \delta}$
denotes the cdf of $U$, and $c$ is the $1- \alpha $ quantile of
$G_{u,v,0}$.  Now $F_{u,\delta}(cuV/v)$ is the probability content
of an origin centered hypersphere of radius $r=\sqrt{cuV/v} $
under translation of an origin centered spherical multivariate
normal distribution by a distance $\delta $. That this is monotone
decreasing in $\delta $ may be taken as obvious by some. It
follows from more general results in Anderson (1955). For a
direct proof start with the fact that
\begin{equation}
F_{u,\delta }(r^{2} )=\delta ^{1-u/2} e^{-\delta ^{2} /2} \int
_{0}^{r}x^{u/2} e^{-x^{2} /2} I_{u/2-1} (\delta x)dx
\end{equation}
where $I_w$ denotes the modified Bessel function of the first kind
and order $w$.  (For example, see Ruben (1962, (3.5)).) The
derivative with respect to \textit{$\delta $} is calculated in
Gilliland (1964) where we find
\begin{equation}
dF_{u,\delta } (r^{2} )/d\delta =-r^{u/2} \delta ^{1-u/2}
e^{-(r^{2}+\delta^{2} )/2} I_{u/2} (r\delta ).
\end{equation}

Since $I_{u/2} (r\delta )>0$ for $r\delta > 0$, we see that
$F_{u,\delta } (r^{2} )$ is a strictly decreasing function of
$\delta $.  Together with (1) and the fact that $\delta=\delta
(\sigma)$ is strictly decreasing in $\sigma $, we see that the
power $\omega $ is a strictly decreasing function of $\sigma$. As
seen, all the necessary monotonicities are easily established if
not very widely known.

It follows that any $100(1 -\gamma)\%$ confidence interval
estimator for $\sigma$ based on the residual sum of squares $V$
transforms to a $100(1-\gamma)\%$ confidence interval estimator
for the power $\omega $. The usual $100(1- \gamma)\%$ CI interval
estimate for $\sigma $ is $a<\sigma<b$ where $a=\sqrt{V/B}$ and
$b=\sqrt{V/A}$ with $A$ and $B$ such that
$F_{v,0}(B)-F_{v,0}(A)=1 -\gamma $ . The typical choices are
$A=F_{v,0}^{-1} (\gamma /2)$ and $B=F_{v,0}^{-1} (1-\gamma /2)$.
The corresponding $100(1-\gamma )\%$ confidence interval estimator
for $\omega $ is
\begin{equation}
1-G_{u,v,\delta (b)} (c)<\omega <1-G_{u,v,\delta (a)} (c).
\end{equation}

Tarasin\'{s}ka (2005, pp. 126-127, Table 1) proposes using
positions $A$ and $B$ to minimize the length of the confidence
interval for the power of the one-sided t-test.  The idea applied
to the power of the F-test would  be to choose $A$ and $B$ to
minimize
\begin{equation}
L = G_{u,v,\delta (\sqrt{V/A})}(c)-G_{u,v,\delta (\sqrt{V/B})} (c)
\end{equation}
subject to the constraint $F_{v,0}(B)-F_{v,0}(A)=1-\gamma $. In
this case, the minimizing $A$ and $B$ depend upon $V$ and the
resulting intervals for power $\omega $ and their corresponding
intervals for $\sigma $ are not shown to have $100(1-\gamma )\%$
coverage probability (See Gilliland and Li (2007)).

Let $Y_{1},...,Y_{n}$ be a random sample from a normal population with mean $%
\mu $ and variance $\sigma ^{2}$. Consider the test for the null
hypothesis $H_{\circ }:\mu =\mu _{\circ }$ against the two-sided
alternative
hypothesis $H_{1}:\mu \neq \mu _{\circ }$ with the level of significance $%
\alpha .$ The null hypothesis $H_{\circ }$ is rejected if
\[
|\bar{Y}-\mu _{\circ }|>\frac{S}{\sqrt{n-1}}\sqrt{c}
\]
in which $\bar{Y}=\dfrac{1}{n}\sum\limits_{i=1}^{n}Y_{i}$, $S^{2}=\dfrac{1}{n%
}\sum\limits_{i=1}^{n}(Y_{i}-\bar{Y})^{2}$ and $c$ are the sample
mean, the sample variance and the $1-\alpha$ quantile of the
central F-distribution $G_{1, n-1, 0}$, respectively.

The power function of the test (Lehmann, 1991)  is defined as:
\begin{equation}
P( |\bar{Y}-\mu _{\circ }|>\frac{S}{\sqrt{n-1}}\sqrt{c})
=1-G_{1,n-1,\delta(\sigma)}(c).
\end{equation}
where $\delta(\sigma)=\dfrac{\sqrt{n}|\mu -\mu _{\circ
}|}{\sigma}$ is non-centrality parameter. By using the invariant
properties of Maximum Likelihood Estimators (MLE), the MLE of
power function in (6) is
\begin{equation}
1-G_{1,n-1,\delta(S)}(c).
\end{equation}

Let (a,b) be any $100(1-\gamma )\%$ confidence interval for
$\sigma $, then by (4) the  confidence interval for power
function is
\begin{equation}
\left\{1-G_{1,n-1,\delta(b)}(c),1-G_{1,n-1,\delta(a)}(c)\right\}.
\end{equation}

The curve of the confidence intervals  in (8) and  the MLE
estimates  in (7) for the power function,  as functions of
$\dfrac{\mu -\mu _{\circ }}{S}$ for $n=10$ and $\alpha =0.05$,
are given in Figure 1.

\begin{figure}[tbp]
\epsfxsize 9 cm \centerline{\epsfbox{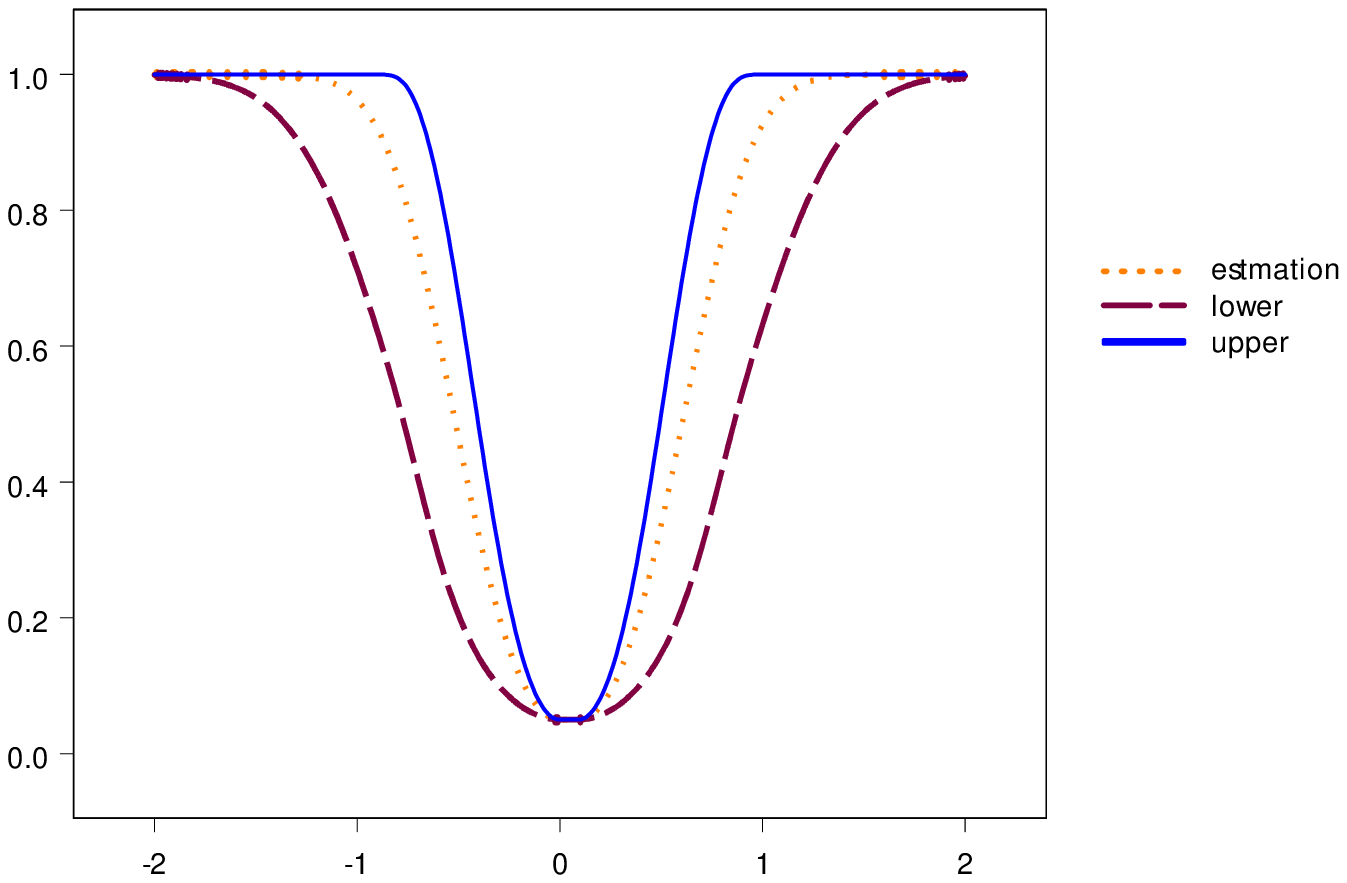}}
\par
\begin{center}
{\small \textbf{Figure 1.} CI curves for the test power and the
estimate of the power as the functions of $\dfrac{\mu -\mu
_{\circ }}{S}
$}%
\end{center}
\end{figure}

\bigskip

 \noindent\textbf{Acknowledgement}: We are grateful to the
editor and referee for their helpful comments and suggestions
which improved the presentation of the result.

\end{document}